# ON THE OPTIMAL DIVIDEND PROBLEM FOR A SPECTRALLY NEGATIVE LÉVY PROCESS

By Florin Avram,[1,2] Zbigniew Palmowski[2,3] and
Martijn R. Pistorius[1]

*Université de Pau, University of Wrocław and King's College London*

In this paper we consider the optimal dividend problem for an insurance company whose risk process evolves as a spectrally negative Lévy process in the absence of dividend payments. The classical dividend problem for an insurance company consists in finding a dividend payment policy that maximizes the total expected discounted dividends. Related is the problem where we impose the restriction that ruin be prevented: the beneficiaries of the dividends must then keep the insurance company solvent by bail-out loans. Drawing on the fluctuation theory of spectrally negative Lévy processes we give an explicit analytical description of the optimal strategy in the set of barrier strategies and the corresponding value function, for either of the problems. Subsequently we investigate when the dividend policy that is optimal among all admissible ones takes the form of a barrier strategy.

**1. Introduction.** In classical collective risk theory (e.g., [11]) the surplus $X = \{X_t, t \geq 0\}$ of an insurance company with initial capital $x$ is described by the Cramér–Lundberg model:

$$(1.1) \qquad X_t = x + \mathrm{d}t - \sum_{k=1}^{N_t} C_k,$$

where $C_k$ are i.i.d. positive random variables representing the claims made, $N = \{N_t, t \geq 0\}$ is an independent Poisson process modeling the times at which the claims occur and $\mathrm{d}t$ represents the premium income up to time

Received January 2006; revised June 2006.
[1]Supported by London Mathematical Society Grant 4416.
[2]Supported by POLONIUM Grant 09158SD.
[3]Supported by KBN 1P03A03128 and NWO 613.000.310.

$t$. Under the assumption that the premium income per unit time d is larger than the average amount claimed, $\lambda \mathbf{E}[C_1]$, the surplus in the Cramér–Lundberg model has positive first moment and has therefore the unrealistic property that it converges to infinity with probability 1. In answer to this objection De Finetti [10] introduced the dividend barrier model, in which all surpluses above a given level are transferred to a beneficiary. In the mathematical finance and actuarial literature there is a good deal of work on dividend barrier models and the problem of finding an optimal policy for paying out dividends. Gerber and Shiu [12] and Jeanblanc and Shiryaev [15] consider the optimal dividend problem in a Brownian setting. Irbäck [14] and Zhou [25] study constant barriers under the model (1.1). Asmussen, Højgaard and Taksar [3] investigate excess-of-loss reinsurance and dividend distribution policies in a diffusion setting. Azcue and Muler [1] follow a viscosity approach to investigate optimal reinsurance and dividend policies in the Cramér–Lundberg model.

A drawback of the dividend barrier model is that under this model the risk process will down-cross the level zero with probability 1. Several ways to combine dividend and ruin considerations are possible; here, we choose one studied in a Brownian motion setting by Harrison and Taylor [13] and Løkka and Zervos [19] involving bail-out loans to prevent ruin, over an infinite horizon.

In this paper we shall approach the dividend problem from the point of view of a general spectrally negative Lévy process. Drawing on the fluctuation theory for spectrally negative Lévy processes, we derive in Sections 3 and 4 expressions for the expectations of the discounted accumulated local time of a reflected and doubly reflected spectrally negative Lévy process, in terms of the scale functions of the Lévy process. Together with known results from the fluctuation theory of spectrally negative Lévy processes and control theory we apply these results in Section 5 to investigate the optimality of barrier dividend strategies for either of the dividend problems. Finally we conclude the paper with some explicit examples in the classical and "bail-out" setting.

**2. Problem setting.** Let $X = \{X_t, t \geq 0\}$ be a Lévy process without positive jumps; that is, $X$ is a stationary stochastic process with independent increments that has right-continuous paths with left limits, only negative jumps and starts at $X_0 = 0$, defined on some filtered probability space $(\Omega, \mathcal{F}, \mathbf{F} = \{\mathcal{F}_t\}_{t \geq 0}, \mathbf{P})$, where $\mathbf{F} = \{\mathcal{F}_t\}_{t \geq 0}$ is a filtration that satisfies the usual conditions of right-continuity and completeness. Denote by $\{\mathbf{P}_x, x \in \mathbf{R}\}$ the family of probability measures corresponding to a translation of $X$ such that $X_0 = x$, where we write $\mathbf{P} = \mathbf{P}_0$. Let $\mathbf{E}_x$ be expectation with respect to $\mathbf{P}_x$. To avoid trivialities, we exclude the case that $X$ has monotone paths. For background on Lévy processes we refer to [24] and [6].

The process $X$ models the risk process of an insurance company or the cash fund of an investment company before dividends are deducted. Let $\pi$ be a dividend strategy consisting of a nondecreasing left-continuous $\mathbf{F}$-adapted process $\pi = \{L_t^\pi, t \geq 0\}$ with $L_0^\pi = 0$, where $L_t^\pi$ represents the cumulative dividends paid out by the company up until time $t$. The risk process with initial capital $x \geq 0$ and controlled by a dividend policy $\pi$ is then given by $U^\pi = \{U_t^\pi, t \geq 0\}$, where

$$(2.1) \qquad U_t^\pi = X_t - L_t^\pi,$$

with $X_0 = x$. Writing $\sigma^\pi = \inf\{t \geq 0 : U_t^\pi < 0\}$ for the time at which ruin occurs, a dividend strategy is called admissible if, at any time before ruin, a lump sum dividend payment is smaller than the size of the available reserves: $L_{t^+}^\pi - L_t^\pi < U_t^\pi$ for $t < \sigma^\pi$. Denoting the set of all admissible strategies by $\Pi$, the expected value discounted at rate $q > 0$ associated to the dividend policy $\pi \in \Pi$ with initial capital $x \geq 0$ is given by

$$v_\pi(x) = \mathbf{E}_x\left[\int_0^{\sigma^\pi} e^{-qt}\,dL_t^\pi\right].$$

The objective of the beneficiaries of the insurance company is to maximize $v_\pi(x)$ over all admissible strategies $\pi$:

$$(2.2) \qquad v_*(x) = \sup_{\pi \in \Pi} v_\pi(x).$$

Consider next the situation where the insurance company is not allowed to go bankrupt and the beneficiary of the dividends is required to inject capital into the insurance company to ensure its risk process stays nonnegative. In this setting a dividend policy $\bar\pi = \{L^{\bar\pi}, R^{\bar\pi}\}$ is a pair of nondecreasing $\mathbf{F}$-adapted processes with $R_0^\pi = L_0^\pi = 0$ such that $R^{\bar\pi} = \{R_t^{\bar\pi}, t \geq 0\}$ is a right-continuous process describing the cumulative amount of injected capital and $L^{\bar\pi} = \{L_t^{\bar\pi}, t \geq 0\}$ is a left-continuous process representing the cumulative amount of paid dividends. Under policy $\bar\pi$ the controlled risk process with initial reserves $x \geq 0$ satisfies

$$V_t^{\bar\pi} = X_t - L_t^{\bar\pi} + R_t^{\bar\pi},$$

where $X_0 = x$. The set of admissible policies $\overline{\Pi}$ consists of those policies for which $V_t^{\bar\pi}$ is nonnegative for $t > 0$ and

$$(2.3) \qquad \int_0^\infty e^{-qt}\,dR_t^{\bar\pi} < \infty, \qquad \mathbf{P}_x\text{-almost surely.}$$

The value associated to the strategy $\bar\pi \in \overline{\Pi}$ starting with capital $x \geq 0$ is then given by

$$\bar v_{\bar\pi}(x) = \mathbf{E}_x\left[\int_0^\infty e^{-qt}\,dL_t^{\bar\pi} - \varphi\int_0^\infty e^{-qt}\,dR_t^{\bar\pi}\right],$$

where $\varphi$ is the cost per unit injected capital, and the associated objective then reads as

$$(2.4) \qquad \bar{v}_*(x) = \sup_{\bar{\pi} \in \overline{\Pi}} \bar{v}_{\bar{\pi}}(x).$$

To ensure that the value function is finite and to avoid degeneracies, we assume that $\mathbf{E}_x[X_1] > -\infty$, $q > 0$ and $\varphi > 1$. To illustrate what happens if $\varphi$ is (close to) 1, we consider the case that $\varphi = 1$ and $X$ is given by (1.1). In this setting, it is no more expensive to pay incoming claims from the reserves or by a bail-out loan, and therefore, as a consequence of the positive discount-factor $q > 0$, it is optimal to pay out all reserves and premiums immediately as dividends and to pay all claims by bail-out loans.

A subclass of possible dividend policies for (2.2), denoted by $\Pi_{\leq C}$, is formed by the set of all strategies $\pi \in \Pi$ under which the controlled risk process $U^\pi$ stays below the constant level $C \geq 0$, $U^\pi(t) \leq C$ for all $t > 0$. An example of an element in $\Pi_{\leq C}$ is a *constant barrier strategy* $\pi_a$ at level $a \leq C$ that corresponds to reducing the risk process $U$ to the level $a$ if $x > a$, by paying out the amount $(x-a)^+$, and subsequently paying out the minimal amount of dividends to keep the risk process below the level $a$. Similarly, in problem (2.4), the double barrier strategy $\bar{\pi}_{0,a}$ with a lower barrier at zero and an upper barrier at level $a$ consists in extracting the required amount of capital to bring the risk process down to the level $a$ and subsequently paying out or in the minimal amount of capital required to keep the risk process between 0 and $a$. In the next section we shall use fluctuation theory of spectrally negative Lévy processes to identify the value functions in problems (2.2) and (2.4) corresponding to the constant barrier strategies $\pi_a$ and $\bar{\pi}_{0,a}$.

**3. Reflected Lévy processes.** We first review some fluctuation theory of spectrally negative Lévy processes and refer the reader for more background to [6, 7, 8, 16, 20, 21] and references therein.

3.1. *Preliminaries.* Since the jumps of a spectrally negative Lévy process $X$ are all nonpositive, the moment generating function $\mathbf{E}[e^{\theta X_t}]$ exists for all $\theta \geq 0$ and is given by $\mathbf{E}[e^{\theta X_t}] = e^{t\psi(\theta)}$ for some function $\psi(\theta)$ that is well defined at least on the positive half-axes where it is strictly convex with the property that $\lim_{\theta \to \infty} \psi(\theta) = +\infty$. Moreover, $\psi$ is strictly increasing on $[\Phi(0), \infty)$, where $\Phi(0)$ is the largest root of $\psi(\theta) = 0$. We shall denote the right-inverse function of $\psi$ by $\Phi : [0, \infty) \to [\Phi(0), \infty)$.

For any $\theta$ for which $\psi(\theta) = \log \mathbf{E}[\exp \theta X_1]$ is finite we denote by $\mathbf{P}^\theta$ an exponential tilting of the measure $\mathbf{P}$ with Radon–Nikodym derivative with respect to $\mathbf{P}$ given by

$$(3.1) \qquad \left.\frac{d\mathbf{P}^\theta}{d\mathbf{P}}\right|_{\mathcal{F}_t} = \exp(\theta X_t - \psi(\theta)t).$$

Under the measure $\mathbf{P}^\theta$ the process $X$ is still a spectrally negative Lévy process with characteristic function $\psi_\theta$ given by

$$(3.2) \qquad \psi_\theta(s) = \psi(s + \theta) - \psi(\theta).$$

Denote by $\sigma$ the Gaussian coefficient and by $\nu$ the Lévy measure of $X$. We recall that if $X$ has bounded variation it takes the form $X_t = \mathrm{d}t - S_t$ for a subordinator $S$ and constant $\mathrm{d} > 0$, also referred to as the infinitesimal drift of $X$. Throughout the paper we assume that $X$ has unbounded variation or $\nu$ is absolutely continuous with respect to the Lebesgue measure:

$$(3.3) \qquad \sigma > 0 \quad \text{or} \quad \int_{-1}^0 x\nu(dx) = \infty \quad \text{or} \quad \nu(dx) \ll dx.$$

3.2. *Scale functions.* For $q \geq 0$, there exists a function $W^{(q)} : [0, \infty) \to [0, \infty)$, called the *q-scale function*, that is continuous and increasing with Laplace transform

$$(3.4) \qquad \int_0^\infty e^{-\theta x} W^{(q)}(y)\, dy = (\psi(\theta) - q)^{-1}, \qquad \theta > \Phi(q).$$

The domain of $W^{(q)}$ is extended to the entire real axis by setting $W^{(q)}(y) = 0$ for $y < 0$. For later use we mention some properties of the function $W^{(q)}$ that have been obtained in the literature. On $(0, \infty)$ the function $y \mapsto W^{(q)}(y)$ is right- and left-differentiable and, as shown in [18], under condition (3.3), it holds that $y \mapsto W^{(q)}(y)$ is continuously differentiable for $y > 0$. The value of the scale function and its derivative in zero can be derived from the Laplace transform (3.4) to be equal to

$$(3.5) \qquad W^{(q)}(0) = 1/\mathrm{d} \quad \text{and} \quad W^{(q)\prime}(0^+) = (q + \nu(-\infty, 0))/\mathrm{d}^2,$$

if $X$ has bounded variation, and $W^{(q)}(0) = W(0) = 0$ if $X$ has unbounded variation (see, e.g., [16], Exercise 8.5 and Lemma 8.3). Moreover, if $\sigma > 0$, it holds that $W^{(q)} \in C^2(0, \infty)$ with $W^{(q)\prime}(0^+) = 2/\sigma^2$; if $X$ has unbounded variation but $\sigma = 0$, it holds that $W^{(q)\prime}(0^+) = \infty$ (see [21], Lemma 4 and [22], Lemma 1).

The function $W^{(q)}$ plays a key role in the solution of the two-sided exit problem as shown by the following classical identity. Letting $T_a^+, T_a^-$ be the entrance times of $X$ into $(a, \infty)$ and $(-\infty, -a)$, respectively:

$$T_a^+ = \inf\{t \geq 0 : X_t > a\}, \qquad T_a^- = \inf\{t \geq 0 : X_t < -a\},$$

and letting $T_{0,a} = T_0^- \wedge T_a^+$ be the first exit time from $[0, a]$, it holds for $y \in [0, a]$ that

$$(3.6) \qquad \mathbf{E}_y[\exp(-q T_{0,a}) \mathbf{1}_{\{T_0^- > T_a^+\}}] = W^{(q)}(y)/W^{(q)}(a),$$

where $\mathbf{1}_A$ is the indicator of the event $A$. Closely related to $W^{(q)}$ is the function $Z^{(q)}$ given by

$$Z^{(q)}(y) = 1 + q\overline{W}^{(q)}(y),$$

where $\overline{W}^{(q)}(y) = \int_0^y W^{(q)}(z)\,dz$ is the antiderivative of $W^{(q)}$. The name $q$-scale function for $W^{(q)}$ and $Z^{(q)}$ is justified as these functions are harmonic for the process $X$ killed upon entering $(-\infty, 0)$, in the sense that

$$\{e^{-q(t\wedge T_0^-)}Z^{(q)}(X_{t\wedge T_0^-}), t \geq 0\} \quad \text{and} \quad \{e^{-q(t\wedge T_0^-)}W^{(q)}(X_{t\wedge T_0^-}), t \geq 0\}$$
(3.7)

are martingales, as shown in [21], Proposition 3. Appealing to this martingale property one can show the following relation between $W^{(q)}$ and its anti-derivative:

LEMMA 1.   *For $y \in [0, a]$ and $a > 0$ it holds that*

$$\overline{W}^{(q)}(y)/\overline{W}^{(q)}(a) \leq W^{(q)}(y)/W^{(q)}(a).$$

PROOF.   Writing $h(y) = \overline{W}^{(q)}(y)/\overline{W}^{(q)}(a) - W^{(q)}(y)/W^{(q)}(a)$ as

$$h(y) = q^{-1}Z^{(q)}(y)/\overline{W}^{(q)}(a) - W^{(q)}(y)/W^{(q)}(a) - q^{-1}/\overline{W}^{(q)}(a)$$

and using the martingale property of $Z^{(q)}$ and $W^{(q)}$ in conjunction with the optional stopping theorem, it follows that $\{e^{-q(t\wedge T_{0,a})}h(X_{t\wedge T_{0,a}}), t \geq 0\}$ can be written as the sum of a martingale and an increasing process and is thus a submartingale. Therefore

$$h(y) \leq \mathbf{E}_y[e^{-q(t\wedge T_{0,a})}h(X_{t\wedge T_{0,a}})] \leq \mathbf{E}_y[e^{-qT_{0,a}}h(X_{T_{0,a}})] = 0,$$

where the last equality follows since $h(y) = 0$ for $y \in (-\infty, 0] \cup \{a\}$.   □

3.3. *Reflection at the supremum.*   Write $I$ and $S$ for the running infimum and supremum of $X$, respectively; that is,

$$I_t = \inf_{0 \leq s \leq t}(X_s \wedge 0) \quad \text{and} \quad S_t = \sup_{0 \leq s \leq t}(X_s \vee 0),$$

where we used the notation $c \vee 0 = \max\{c, 0\}$ and $c \wedge 0 = \min\{c, 0\}$. By $Y = X - I$ and $\widehat{Y} = S - X$ we denote the Lévy process $X$ reflected at its past infimum $I$ and at its past supremum $S$, respectively. Denoting by $\eta(q)$ an independent random variable with parameter $q$, it follows, by duality and the Wiener–Hopf factorization of $X$ (e.g., [6], page 45 and pages 188–192, resp.), that

(3.8)                $S_{\eta(q)} \sim Y_{\eta(q)} \sim \exp(\Phi(q)).$

Further, it was shown in [4] and [21] that the Laplace transform of the entrance time $\tau_a$ of the reflected process $Y$ into $(a, \infty)$ [resp. the entrance time $\hat{\tau}_a$ of $\hat{Y}$ into $(a, \infty)$] can be expressed in terms of the functions $Z^{(q)}$ and $W^{(q)}$ as follows:

$$(3.9) \qquad \mathbf{E}_y[e^{-q\tau_a}] = \frac{Z^{(q)}(y)}{Z^{(q)}(a)},$$

$$(3.10) \qquad \mathbf{E}_{-y}[e^{-q\hat{\tau}_a}] = Z^{(q)}(a-y) - qW^{(q)}(a-y)\frac{W^{(q)}(a)}{W^{(q)\prime}(a)},$$

where $y \in [0, a]$ and where we note that under $\mathbf{P}_y$ [$\mathbf{P}_{-y}$] it holds that $Y_0 = y$ [$\hat{Y}_0 = y$]. The identity (3.10) together with the strong Markov property implies the martingale property of

$$(3.11) \qquad e^{-q(t \wedge \hat{\tau}_a)}\left\{Z^{(q)}(a - \hat{Y}_{t \wedge \hat{\tau}_a}) - qW^{(q)}(a - \hat{Y}_{t \wedge \hat{\tau}_a})\frac{W^{(q)}(a)}{W^{(q)\prime}(a)}\right\}.$$

Denote by $\pi_a = \{L_t^a, t \leq \sigma_a\}$ the constant barrier strategy at level $a$ and let $U^a = U^{\pi_a}$ be the corresponding risk process. If $U_0^a \in [0, a]$, the strategy $\pi_a$ corresponds to a reflection of the process $X - a$ at its supremum: for $t \leq \sigma_a$ process $L_t^a$ can be explicitly represented by

$$L_t^a = \sup_{s \leq t}[X_s - a] \vee 0.$$

Note that $L_t^a$ is a Markov local time of $U^a$ at $a$, that is, $L^a$ is increasing, continuous and adapted such that the support of the Stieltjes measure $dL_t^a$ is contained in the closure of the set $\{t : U_t^a = a\}$ (see, e.g. [6], Chapter IV, for background on local times). In the case that $U_0 = x > a$, $L_t^a$ has a jump at $t = 0$ of size $x - a$ to bring $U^a$ to the level $a$ and a similar structure afterward: $L_t^a = (x-a)\mathbf{1}_{\{t>0\}} + \sup_{s \leq t}[X_s - x] \vee 0$.

The following result concerns the value function associated to the dividend policy $\pi_a$:

PROPOSITION 1. *Let $a > 0$. For $x \in [0, a]$ it holds that*

$$(3.12) \qquad \mathbf{E}_x\left[\int_0^{\sigma_a} e^{-qt} \, dL_t^a\right] = \mathbf{E}_{x-a}\left[\int_0^{\hat{\tau}_a} e^{-qt} \, dS_t\right] = \frac{W^{(q)}(x)}{W^{(q)\prime}(a)},$$

*where $\sigma_a = \sigma^{\pi_a} = \inf\{t \geq 0 : U_t^a < 0\}$ is the ruin time.*

PROOF. By spatial homogeneity of the Lévy process $X$, it follows that the ensemble $\{U_t^a, L_t^a, t \leq \sigma_a; U_0 = x\}$ has the same law as $\{a - \hat{Y}_t, S_t, t \leq \hat{\tau}_a; \hat{Y}_0 = a - x\}$. Noting that $\hat{Y}_0 = a - x$ precisely if $X_0 = x - a$ (since then $S_0 = 0$), the first equality of (3.12) is seen to hold true. Using excursion theory it was shown in the proof of [4], Theorem 1, that

$$(3.13) \qquad \mathbf{E}_0\left[\int_0^{\hat{\tau}_a} e^{-qt}\, dS_t\right] = \frac{W^{(q)}(a)}{W^{(q)\prime}(a)}.$$

Applying the strong Markov property of $\widehat{Y}$ at $\hat{\tau}_0 = \inf\{t \geq 0 : \widehat{Y}_t = 0\}$ and using that $\{\widehat{Y}_t, t \leq \hat{\tau}_0\}$ is in law equal to $\{-X_t, t \leq T_0^+, X_0 = -\widehat{Y}_0 = x - a\}$, we find that

$$\mathbf{E}_{x-a}\left[\int_0^{\hat{\tau}_a} e^{-qt}\, dS_t\right] = \mathbf{E}_{x-a}[e^{-q\hat{\tau}_0}\mathbf{1}_{\{\hat{\tau}_0 < \hat{\tau}_a\}}]\mathbf{E}_0\left[\int_0^{\hat{\tau}_a} e^{-qt}\, dS_t\right]$$

$$= \mathbf{E}_{x-a}[e^{-qT_0^+}\mathbf{1}_{\{T_0^+ < T_a^-\}}]\mathbf{E}_0\left[\int_0^{\hat{\tau}_a} e^{-qt}\, dS_t\right].$$

Inserting the identities (3.13) and (3.6) into this equation completes the proof. $\quad\square$

Let us complement the previous result by considering what happens in the case that the barrier is taken to be 0. If $X$ has unbounded variation, 0 is regular for $(-\infty, 0)$ so that $U^0$ immediately enters the negative half-axis and $\mathbf{P}_0(\sigma_0 = 0) = 1$, and the right-hand side of (3.12) is zero (if $x = a = 0$). If $\nu(-\infty, 0)$ is finite, $U^0$ enters $(-\infty, 0)$ when the first jump occurs so that $\sigma_0$ is exponential with mean $\nu(-\infty, 0)^{-1}$ and

$$(3.14) \qquad \mathbf{E}_0\left[\int_0^{\sigma_0} e^{-qt}\, dL_t^0\right] = \mathrm{d}\mathbf{E}_0\left[\int_0^{\sigma_0} e^{-qt}\, dt\right] = \frac{\mathrm{d}}{q + \nu(-\infty, 0)}.$$

If $\nu$ is infinite but $X$ has bounded variation, the validity of (3.14) follows by approximation. Combining these observations with (3.5), we note that (3.12) remains valid for $x = a = 0$ if $W^{(q)\prime}(a)$ for $a = 0$ is understood to be $W^{(q)\prime}(0^+)$.

In view of (3.8), (3.13) and since $a \mapsto \hat{\tau}_a$ is nondecreasing with $\lim_{a\to\infty}\hat{\tau}_a = +\infty$ a.s., we note for later reference that $W^{(q)}/W^{(q)\prime}$ is an increasing function on $(0, \infty)$ with limit

$$(3.15) \qquad \lim_{a\to\infty}\frac{W^{(q)}(a)}{W^{(q)\prime}(a)} = \mathbf{E}_0[S_{\eta(q)}] = \frac{1}{\Phi(q)}.$$

3.4. *Martingales and overshoot.* In the sequel we shall need the following identities of expected discounted overshoots and related martingales in terms of the anti-derivative $\overline{Z}^{(q)}(y)$ of $Z^{(q)}(y)$ which is for $y \in \mathbf{R}$ defined by

$$\overline{Z}^{(q)}(y) = \int_0^y Z^{(q)}(z)\, dz = y + q\int_0^y \int_0^z W^{(q)}(w)\, dw\, dz.$$

Note that $\overline{Z}^{(q)}(y) = y$ for $y < 0$, since we set $W^{(q)}(y) = 0$ for $y < 0$.

PROPOSITION 2. *If $\psi'(0^+) > -\infty$, then*

$$e^{-q(t \wedge T_0^-)}\{\overline{Z}^{(q)}(X_{t \wedge T_0^-}) + \psi'(0^+)/q\}$$

*and*

$$e^{-q(t \wedge \hat{\tau}_a)}\{\overline{Z}^{(q)}(a - \widehat{Y}_{t \wedge \hat{\tau}_a}) + \psi'(0^+)/q - W^{(q)}(a - \widehat{Y}_{t \wedge \hat{\tau}_a})Z^{(q)}(a)/W^{(q)\prime}(a)\}$$

*are martingales. In particular, it holds that for $y \in [0, a]$ and $x \geq 0$,*

$$(3.16) \quad \mathbf{E}_{y-a}[e^{-q\hat{\tau}_a}(a - \widehat{Y}_{\hat{\tau}_a})] = \overline{Z}^{(q)}(y) - \psi'(0^+)\overline{W}^{(q)}(y) - CW^{(q)}(y),$$

$$(3.17) \quad \mathbf{E}_x[e^{-qT_0^-}X_{T_0^-}] = \overline{Z}^{(q)}(x) - \psi'(0^+)\overline{W}^{(q)}(x) - DW^{(q)}(x),$$

*where $D = [q - \psi'(0^+)\Phi(q)]/\Phi(q)^2$ and $C = [Z^{(q)}(a) - \psi'(0^+)W^{(q)}(a)]/W^{(q)\prime}(a)$.*

PROOF. We first show the validity of the identities (3.16) and (3.17). Writing $W_v^{(q)}$ and $Z_v^{(q)}$ for the ("tilted") $q$-scale functions of $X$ under $\mathbf{P}^v$ we read off from [4], Theorem 1 and [17], Theorem 4, that for $\kappa := q - \psi(v) \geq 0$, $x \geq 0$ and $y \leq a$ it holds that

$$(3.18) \quad \mathbf{E}_{y-a}[e^{-q\hat{\tau}_a - v(\widehat{Y}_{\hat{\tau}_a} - a)}] = e^{vy}[Z_v^{(\kappa)}(y) - C_v W_v^{(\kappa)}(y)],$$

$$(3.19) \quad \mathbf{E}_x[e^{-qT_0^- + vX_{T_0^-}}] = e^{vx}[Z_v^{(\kappa)}(x) - D_v W_v^{(\kappa)}(x)],$$

where $D_v = \kappa/[\Phi(q) - v]$ and $C_v = [\kappa W_v^{(\kappa)}(a) + vZ_v^{(\kappa)}(a)]/[W_v^{(\kappa)\prime}(a) + vW_v^{(\kappa)}(a)]$. The "tilted" scale functions can be linked to nontilted scale functions via the relation $e^{vy}W_v^{(q-\psi(v))}(y) = W^{(q)}(y)$ from [4], Remark 4. This relation implies that $e^{vy}[W_v^{(\kappa)\prime}(y) + vW_v^{(\kappa)}(y)] = W^{(q)\prime}(y)$ and

$$Z_v^{(\kappa)}(y) = 1 + \kappa \int_0^y e^{-vz}W^{(q)}(z)\,dz.$$

In view of these relations it is a matter of algebra to verify that the right-derivatives with respect to $v$ in $v = 0$ of $D_v$, $C_v$ and $\ell(v) := e^{vy}Z_v^{(q-\psi(v))}(y)$ are respectively equal to the constants $D$ and $C$ given in the statement of the proposition and

$$\ell'(0^+) = \overline{Z}^{(q)}(y) - \psi'(0^+)\overline{W}^{(q)}(y).$$

Differentiating (3.18) and (3.19) and inserting the derived results we arrive at the equations (3.16) and (3.17).

Write now $h_1, h_2$ for the right-hand sides of (3.16) and (3.17), respectively. From the overshoot identities (3.16) and (3.17) and the definition for $y < 0$ of $W^{(q)}(y)$, $\overline{Z}^{(q)}(y)$ and $\overline{W}^{(q)}(y)$, it is straightforward to verify that

$$\mathbf{E}_x[e^{-qT_0^-}h_1(X_{T_0^-})] = h_1(x), \qquad \mathbf{E}_{y-a}[e^{-q\hat{\tau}_a}h_2(a - \widehat{Y}_{\hat{\tau}_a})] = h_2(y).$$

The strong Markov property then implies that for $t \geq 0$,

$$\mathbf{E}_x[e^{-qT_0^-}h_1(X_{T_0^-})|\mathcal{F}_t] = e^{-q(t \wedge T_0^-)}h_1(X_{t \wedge T_0^-}),$$

$$\mathbf{E}_{y-a}[e^{-q\hat{\tau}_a}h_2(a - \hat{Y}_{\hat{\tau}_a})|\mathcal{F}_t] = e^{-q(t \wedge \hat{\tau}_a)}h_2(a - \hat{Y}_{t \wedge \hat{\tau}_a})$$

and, in view of (3.7) and (3.11), the stated martingale properties follow. $\square$

**4. Doubly reflected Lévy processes.** Now we turn to the computation of the value function corresponding to the constant barrier strategy $\bar{\pi}_{0,a} = \{L_t^a, R_t^0, t \geq 0\}$ that consists of imposing "reflecting" barriers $L^a$ and $R^0$ at $a$ and $0$, respectively. When the initial capital $X_0 = x \in [0, a]$ the risk process $V_t^a := V_t^{\bar{\pi}_{0,a}}$ is a doubly reflected spectrally negative Lévy process. Informally, this process moves as a Lévy process while it is inside $[0, a]$ but each time it attempts to down-cross $0$ or up-cross $a$ it is "regulated" to keep it inside the interval $[0, a]$. In [20] a pathwise construction of a doubly reflected Lévy process was given, showing that $V^a$ is a strong Markov process. See also [2], XIV.3, for a discussion of processes with two reflecting barriers in the context of queueing models. It was shown ([20], Theorem 1) that a version of the $q$-potential measure $\widetilde{U}^q(x, dy) = \int_0^\infty e^{-qt}\mathbf{P}_x(V_t^a \in dy)$ of $V^a$ is given by $\widetilde{U}^q(x, dy) = \tilde{u}^q(x, 0)\delta_0(dy) + \tilde{u}^q(x, y)\,dy$ where $\delta_0$ is the pointmass in zero, $\tilde{u}^q(x, 0) = Z^{(q)}(a - x)W^{(q)}(0)/(qW^{(q)}(a))$ and

$$\tilde{u}^q(x, y) = \frac{Z^{(q)}(a - x)W^{(q)\prime}(y)}{qW^{(q)}(a)} - W^{(q)}(y - x),$$

(4.1)
$$x, y \in [0, a], y \neq 0.$$

For $t \geq 0$, $V_t^a$ can be expressed in terms of $X, L^a$ and $R^0$ as

(4.2)
$$V_t^a = X_t - L_t^a + R_t^0$$

for some increasing adapted processes $L^a$ and $R^0$ such that the supports of the Stieltjes measures $dL_t^a$ and $dR_t^0$ are contained in the closures of the sets $\{t : V_t^a = a\}$ and $\{t : V_t^a = 0\}$, respectively. For completeness we extend the construction in [20] to a simultaneous construction of the processes $L^a$, $R^0$ and $V^a$, when $X_0 \in [0, a]$:

0. Set $\sigma = T_{0,a}$. For $t < \sigma$, set $L_t^a = R_t^0 = 0$ and $V_t^a = X_t$.

   If $X_\sigma \leq 0$ set $\xi := X_\sigma$ and go to step 2; else set $L_\sigma^a = 0$ and $V_\sigma^a = a$ and go to step 1.

1. Set $Z_t = X_t - X_\sigma$. For $\sigma < t < \sigma' := \inf\{u \geq \sigma : Z_u < -a\}$, set

   $$L_t^a = L_\sigma^a + \sup_{\sigma \leq s \leq t}[Z_s \vee 0], \qquad V_t^a = a + Z_t - (L_t^a - L_\sigma^a)$$

   and let $R_t^0 = R_\sigma^0$. Set $\sigma := \sigma'$ and $\xi := X_{\sigma'} - X_\sigma + a$ and go to step 2.

2. Set $Z_t = X_t - X_\sigma$. For $\sigma \le t \le \sigma'' := \inf\{u \ge \sigma : Z_u = a\}$, set

$$R_t^0 = R_{\sigma-}^0 - \xi - \inf_{\sigma \le s \le t}[Z_s \wedge 0], \qquad V_t^a = Z_t + R_t^0 - R_\sigma^0$$

and let $L_t^a = L_\sigma^a$. Set $\sigma := \sigma''$ and go to step 1.

It can be verified by induction that the process $V$ constructed in this way satisfies $V_t \in [0, a]$ and $L^a$ and $R^0$ are processes with the required properties such that (4.2) holds.

REMARK. If the initial capital $x > a$, then the above construction can be easily adapted: in step 0 set $L_0^a = 0$, $V_0^a = x$ and $L_{0^+}^a = x - a$, $V_{0^+}^a = a$ and in step 1 set $\sigma = 0$ and replace $L_0^a$ by $L_{0^+}^a$ and repeat the rest of the construction.

In the next result, the expectations of the Laplace–Stieltjes transforms of $L^a$ and $R^0$ are identified:

THEOREM 1. *Let $a > 0$. For $x \in [0, a]$ and $q > 0$ it holds that*

$$(4.3) \qquad \mathbf{E}_x\left[\int_0^\infty e^{-qt}\, dL_t^a\right] = Z^{(q)}(x)/[qW^{(q)}(a)],$$

$$(4.4) \qquad \mathbf{E}_x\left[\int_0^\infty e^{-qt}\, dR_t^0\right] = -\overline{Z}^{(q)}(x) - \frac{\psi'(0^+)}{q} + \frac{Z^{(q)}(a)}{qW^{(q)}(a)}Z^{(q)}(x),$$

*where the expression in (4.4) is understood to be $+\infty$ if $\psi'(0^+) = -\infty$.*

REMARK. If $X$ has bounded variation we can also consider the strategy of immediately paying out all dividends and paying all incoming claims with bail-out loans—this corresponds to keeping the risk process constant equal to zero. Denoting the "reflecting barriers" corresponding to this case by $L^0$ and $R^0$ respectively, one can directly verify that

$$(4.5) \quad \mathbf{E}_0\left[\int_0^\infty e^{-qt}\, dL_t^a\right] = \mathrm{d}/q, \qquad \mathbf{E}_0\left[\int_0^\infty e^{-qt}\, dR_t^0\right] = (\mathrm{d} - \psi'(0^+))/q.$$

PROOF OF THEOREM 1. We first prove (4.3). Denote by $f(u)$ its left-hand side and write $\tau_0' = \inf\{t \ge 0 : V_t^a = b\}$ for the first hitting time of $\{b\}$. We shall derive a recursion for $f(x)$ by considering one cycle of the process $V^a$. More specifically, applying the strong Markov property of $V^a$ at $\tau_0'$ we find that

$$(4.6) \qquad f(x) = \mathbf{E}_x\left[\int_0^{\tau_0'} e^{-qt}\, dL_t^a\right] + \mathbf{E}_x[e^{-q\tau_0'}]f(0).$$

Since $\{V_t^a, t < \tau_0', V_0^a = x\}$ has the same law as $\{a - \widehat{Y}_t, t < \hat{\tau}_a, \widehat{Y}_0 = a - x\}$, the first term and first factor in the second term in (4.6) are equal to (3.12) and (3.10) (with $y = a - x$), respectively. By the fact that $L^a$ does not increase until $V^a$ reaches the level $a$, we find by the strong Markov property that $f(0) = \mathbf{E}_0[e^{-q\tau'_a}]f(a)$, where $\mathbf{E}_0[e^{-q\tau'_a}] = Z^{(q)}(a)^{-1}$ in view of (3.9) and the fact that $\{V^a_t, t \leq \tau'_a, V^a_0 = x\}$ has the same law as $\{Y_t, t \leq \tau_a, Y_0 = x\}$. Inserting all the three formulas into (4.6) results in the equation

$$(4.7) \qquad f(x) = \frac{W^{(q)}(x)}{W^{(q)\prime}(a)} + f(a)\left(\frac{Z^{(q)}(x)}{Z^{(q)}(a)} - q\frac{W^{(q)}(x)W^{(q)}(a)}{Z^{(q)}(a)W^{(q)\prime}(a)}\right).$$

As this relation remains valid for $x = a$, we are led to a recursion for $f(a)$ the solution of which reads as $f(a) = Z^{(q)}(a)/[qW^{(q)}(a)]$. Inserting $f(a)$ back in (4.7) finishes the proof of (4.3).

Now we turn to the calculation of the expectation (4.4). Writing $g(x)$ for the left-hand side of (4.4) and applying the strong Markov property of $V^a$ at $\tau'_a$ shows that

$$(4.8) \qquad g(x) = \mathbf{E}_x\left[\int_0^{\tau'_a} e^{-qt}\,dR^0_t\right] + \mathbf{E}_x[e^{-q\tau'_a}]g(a)$$

with $g(a) = \mathbf{E}_a[e^{-q\tau'_0}\Delta R^0_{\tau'_0}] + \mathbf{E}_a[e^{-q\tau'_0}]g(0)$, where $\Delta R^0_{\tau'_0} = R^0_{\tau'_0} - R^0_{\tau'_0-}$ denotes the jump of $R$ at $\tau'_0$. Appealing to the fact that $\{V^a_t, t < \tau'_0, V^a_0 = x\}$ and $\{V^a_t, t \leq \tau'_a, V^a_0 = x\}$ have the same distribution as $\{a - \widehat{Y}_t, t < \hat{\tau}_a, \widehat{Y}_0 = a - x\}$ and $\{Y_t, t \leq \tau_a, Y_0 = x\}$, respectively, the Laplace transforms of $\tau'_a, \tau'_0$ and the expectation involving $\Delta R^0_{\tau'_0}$ can be identified by (3.9) (with $y = x$), (3.10) (with $y = 0$) and (3.16) respectively. The rest of the proof is devoted to the computation of the first term on the right-hand side of (4.8). Invoking the strong Markov property shows that

$$
\begin{aligned}
-\mathbf{E}_x\left[\int_0^{\tau'_a} e^{-qt}\,dR^0_t\right] &= \mathbf{E}_x\left[\int_0^{\tau_a} e^{-qt}\,dI_t\right] \\
&= \mathbf{E}_x\left[\int_0^{\infty} e^{-qt}\,dI_t\right] - \mathbf{E}_x[e^{-q\tau_a}]\mathbf{E}_a\left[\int_0^{\infty} e^{-qt}\,dI_t\right] \\
&= k(x) - \frac{Z^{(q)}(x)}{Z^{(q)}(a)}k(a),
\end{aligned}
$$

where $k(x) = \mathbf{E}_x[I_{\eta(q)}]$ satisfies

$$
\begin{aligned}
\mathbf{E}_x[I_{\eta(q)}] &= \mathbf{E}_x[e^{-qT^-_0}X_{T^-_0}] + \mathbf{E}_x[e^{-qT^-_0}]\mathbf{E}_0[I_{\eta(q)}] \\
&= \overline{Z}^{(q)}(x) - \Phi(q)^{-1}Z^{(q)}(x) + \psi'(0^+)/q.
\end{aligned}
$$

In the last line we inserted the identities (3.17) and (3.19) (with $v = 0$). Further, we used that $\mathbf{E}_0[I_{\eta(q)}] = \mathbf{E}_0[X_{\eta(q)}] - \mathbf{E}_0[(X - I)_{\eta(q)}]$ where $\mathbf{E}_0[X_{\eta(q)}] = \psi'(0^+)/q$ (from the definition of $\psi$) and $\mathbf{E}_0[(X - I)_{\eta(q)}] = 1/\Phi(q)$ from (3.8).

Inserting the found identities into (4.8) and taking $x$ to be zero in (4.8) yields a recursion for $g(0)$, which can be solved explicitly in terms of the scale functions. After some algebra one arrives at

$$g(0) = -\frac{\psi'(0^+)}{q} + \frac{Z^{(q)}(a)}{qW^{(q)}(a)}.$$

Substituting this expression back into (4.8) results in (4.4). $\quad\square$

## 5. Optimal dividend strategies.

When solving the optimal dividend problems our method draws on classical optimal control literature such as, for example, Jeanblanc and Shiryaev [15] and Harrison and Taylor [13] who deal with the classical dividend problem and a storage system in a Brownian motion setting, respectively. In these papers it was shown that if the state process follows a Brownian motion with drift, the optimal strategy takes the form of a barrier strategy. In view of the fact that our state process is still a Markov process, we expect that barrier strategies play an important role in the solution of the problems (2.2) and (2.4). In this section we shall investigate the optimality of barrier strategies in the classical dividend problem (2.2) and the bail-out problem (2.4).

5.1. *Classical dividend problem.* From Proposition 1 we read off that the value functions corresponding to barrier strategies $\pi_a$ at the levels $a > 0$ are given by

$$(5.1) \qquad v_a(x) = v_{\pi_a}(x) = \begin{cases} \dfrac{W^{(q)}(x)}{W^{(q)\prime}(a)}, & 0 \le x \le a, \\ x - a + \dfrac{W^{(q)}(a)}{W^{(q)\prime}(a)}, & x > a, \end{cases}$$

and the strategy of taking out all dividends immediately has value $v_0(x) = x + W^{(q)}(0)/W^{(q)\prime}(0^+)$. To complete the description of the candidate optimal barrier solution we specify the level $c^*$ of the barrier as

$$(5.2) \qquad c^* = \inf\{a > 0 : W^{(q)\prime}(a) \le W^{(q)\prime}(x) \text{ for all } x\},$$

where $\inf \varnothing = \infty$. Note that, if $W^{(q)}$ is twice continuously differentiable on $(0, \infty)$ (which is in general not the case) and $c^* > 0$, then $c^*$ satisfies

$$(5.3) \qquad W^{(q)\prime\prime}(c^*) = 0,$$

so that in that case the optimal level $c^*$ is such that the value function $v_{c^*}$ is $C^2$ on $(0, \infty)$. Recalling that $W^{(q)\prime}(0^+)$ is infinite if $X$ has no Brownian component and the mass of its Lévy measure $\nu$ is infinite, we infer from the definition (5.2) of $c^*$ that in this case $c^* > 0$ irrespective of the sign of the drift $\mathbf{E}_x[X_1]$. In comparison, if $X$ is a Brownian motion with drift $\mu$, $c^*$ is positive or zero according to whether the drift $\mu$ is positive or nonpositive. (See also Section 6 for other specific examples.)

Denote by $\Gamma$ the extended generator of the process $X$, which acts on $C^2$ functions $f$ with compact support as

$$\Gamma f(x) = \frac{\sigma^2}{2} f''(x) + cf'(x) + \int_{-\infty}^{0} [f(x+y) - f(x) - f'(x)y \mathbf{1}_{\{|y|<1\}}]\nu(dy),$$

where $\nu$ is the Lévy measure of $X$ and $\sigma^2$ denotes the Gaussian coefficient and $c = \mathrm{d} + \int_{-1}^{0} y\nu(dy)$ if the jump-part has bounded variation; see [9], Theorem 7.14 and [24], Chapter 6, Theorem 31.5. Note that by the properties of $W^{(q)}$ given in Section 3.2, it follows that $v_{c^*}$ is $C^2$ on $(0,\infty)$ if $\sigma > 0$ and is $C^1$ on $(0,\infty)$ if $X$ has bounded variation. The following result concerns optimality of the barrier strategy $\pi_{c^*}$ for the classical dividend problem.

THEOREM 2. *Assume that $\sigma > 0$ or that $X$ has bounded variation or, otherwise, suppose that $v_{c^*} \in C^2(0,\infty)$. If $q > 0$, then $c^* < \infty$ and the following hold true:*

  (i) *$\pi_c^*$ is the optimal strategy in the set $\Pi_{\leq c^*}$ and $v_{c^*} = \sup_{\pi \in \Pi_{\leq c^*}} v_\pi$.*

  (ii) *If $(\Gamma v_{c^*} - qv_{c^*})(x) \leq 0$ for $x > c^*$, the value function and optimal strategy of (2.2) are given by $v_* = v_{c^*}$ and $\pi_* = \pi_{c^*}$, respectively.*

REMARK. If the condition $(\Gamma v_{c^*} - qv_{c^*})(x) \leq 0$ is not satisfied for all $x \geq c^*$, but if $c^* > 0$ and one can construct a function $v$ on $[0,\infty)$ that satisfies the Hamilton–Jacobi–Bellman (HJB) equation (5.8) (see for a precise statement Proposition 5 below), the strategy $\pi_{c^*}$ is optimal for "small" initial reserves, that is, it is optimal to apply the barrier strategy $\pi_{c^*} \circ \theta_t$ whenever $U_t \in [0,c^*]$ (where $\theta$ denotes the shift operator) and it holds that $v(x) = v_*(x) = v_{c^*}(x)$ for $x \in [0,c^*]$. This observation agrees with the description of the optimal value function in the setting of the Cramér–Lundberg model, obtained in [1], Section 9, using viscosity methods. Azcue and Muler [1] also constructed an example with $c^* = 0$ and $(\Gamma v_{c^*} - qv_{c^*})(x) > 0$ for some $x > 0$ where the optimal strategy does not take the form of a barrier strategy.

5.2. *Dividends and bail-out.* In the "bail-out" setting and under the assumption that $\psi'(0^+) > -\infty$, we read off from Theorem 1 that the value function corresponding to the strategy $\bar{\pi}_{0,a}$ of putting reflecting barriers at the levels 0 and $a > 0$ is given by $\bar{v}_{\bar{\pi}_{0,a}} = \bar{v}_a$ where

$$(5.4) \quad \bar{v}_a(x) = \begin{cases} \varphi(\overline{Z}^{(q)}(x) + \psi'(0^+)/q) + Z^{(q)}(x)\left[\dfrac{1 - \varphi Z^{(q)}(a)}{qW^{(q)}(a)}\right], \\ \hspace{6cm} 0 \leq x \leq a, \\ x - a + \bar{v}_a(a), \hspace{3.5cm} x > a. \end{cases}$$

In particular, if $X$ is a Lévy process of bounded variation with drift d,

$$(5.5) \qquad \bar{v}_0(x) = x + [\varphi\psi'(0^+) + (1-\varphi)\mathrm{d}]/q$$

is the value function corresponding to keeping the risk process identically equal to zero. The barrier level is specified as

$$(5.6) \quad d^* = \inf\{a > 0 : G(a) := [\varphi Z^{(q)}(a) - 1]W^{(q)\prime}(a) - \varphi q W^{(q)}(a)^2 \le 0\}.$$

Below, in Lemma 2, we shall show that if $\nu(-\infty, 0) \le \frac{q}{\varphi-1}$ and there is no Brownian component, then $d^* = 0$; else $d^* > 0$.

The constructed solution $\bar{v}_{d^*}$ can be identified as the value function of the optimal dividend problem (2.4). As a consequence, in the bail-out setting the optimal strategy takes the form of a barrier strategy for any initial capital:

THEOREM 3.    *Let $q > 0$ and suppose that $\psi'(0^+) < \infty$. Then $d^* < \infty$ and the value function and optimal strategy of (2.4) are given by $\bar{v}_*(x) = \bar{v}_{d^*}(x)$ and $\bar{\pi}_* = \bar{\pi}_{0,d^*}$, respectively.*

5.3. *Optimal barrier strategies.* As a first step in proving Theorems 2 and 3 we show optimality of $\pi_{c^*}$ and $\bar{\pi}_{0,d^*}$ across the respective set of barrier strategies:

PROPOSITION 3.    *Let $q > 0$.*

(i)  *It holds that $c^* < \infty$ and $\pi_{c^*}$ is an optimal barrier strategy if $U_0 \in [0, c^*]$, that is,*

$$v_a(x) \le v_{c^*}(x), \qquad x \in [0, c^*], a \ge 0.$$

(ii)  *Suppose that $\psi'(0^+) < \infty$. It holds that $d^* < \infty$ and $\bar{\pi}_{0,d^*}$ is the optimal barrier strategy for any initial capital, that is,*

$$\bar{v}_a(x) \le \bar{v}_{d^*}(x), \qquad x, a \ge 0.$$

To prove Proposition 3 we use the following facts regarding $c^*$ and $d^*$:

LEMMA 2.    *Suppose that $q > 0$.*

(i)  *It holds that $c^* < \infty$.*
(ii)  *If $\nu(-\infty, 0) \le q/(\varphi - 1)$ and $\sigma = 0$, then $d^* = 0$; else $d^* > 0$.*

PROOF.    (i) Recall that $W^{(q)\prime}(y)$ is nonnegative and continuous for $y > 0$ and increases to $\infty$ as $y \to \infty$. Therefore it holds that either $W^{(q)\prime}(y)$ attains its finite minimum at some $y \in (0, \infty)$ or $W^{(q)\prime}(0^+) \le W^{(q)\prime}(y)$ for all $y \in (0, \infty)$.

(ii) Write $H(a) = \mathbf{E}_0[e^{-q\hat{\tau}_a}]$ and recall that $H(a)$ is given by (3.10) [with $y = 0$]. As $W^{(q)}(a)$ and $W^{(q)\prime}(a)$ are strictly positive, we see that $G(a) = 0$, with $G$ defined in (5.6), can be rewritten as $F(a) = 0$ where

(5.7)           $$F(a) := [\varphi H(a) - 1]W^{(q)\prime}(a)/[qW^{(q)}(a)^2].$$

Since $a \mapsto \hat{\tau}_a$ is monotonically increasing with $\lim_{a\to\infty} \hat{\tau}_a = \infty$ almost surely, it follows that $H(a)$ monotonically decreases to zero as $a \to \infty$. Therefore $F(a) \leq 0$ for all $a > 0$ if $F(0^+) \leq 0$. Further, as $F$ is continuous, it follows as a consequence of the intermediate value theorem that $F(a) = 0$ has a root in $(0, \infty)$ if $F(0^+) \in (0, \infty]$. In view of the fact that both $W^{(q)}(0^+) > 0$ and $W^{(q)\prime}(0^+) < \infty$ hold true precisely if $X$ is a compound Poisson process, we see that $F(0^+) \leq 0$ if and only if both $\sigma = 0$ and $\nu(-\infty, 0) \leq q/(\varphi - 1)$ are satisfied. The statement (ii) follows.   □

PROOF OF PROPOSITION 3.   (i) From Lemma 2(i) we have that $c^* < \infty$. Combining the definition of $c^*$ and the following estimate for all $x \geq b$,

$$x - b + \frac{W^{(q)}(b)}{W^{(q)\prime}(b)} \leq x - b + \frac{W^{(q)}(b)}{W^{(q)\prime}(c^*)} \leq \frac{W^{(q)}(x)}{W^{(q)\prime}(c^*)},$$

the assertion follows in view of the definition of $v_a$.

(ii) It is straightforward to verify that, for any $x > 0$, the derivative of $a \mapsto \bar{v}_a(x)$ in $a > x$ is equal to $F(a) \times [Z^{(q)}(x)]$ and in $0 < a < x$ is equal to $F(a) \times [Z^{(q)}(a)]$ respectively, where $F(a)$ is given in (5.7). In particular we note that, for any $x > 0$, $a \mapsto \bar{v}_a(x)$ is $C^1$ on $(0, \infty)$. In view of the arguments in the proof of Lemma 2 and the definition of $d^*$ we see that $F(a) \leq 0$ for $a > d^*$, and if $d^* > 0$, $F(d^*) = 0$ and $F(a) > 0$ for $0 < a < d^*$. Therefore $a \mapsto \bar{v}_a(x)$ attains its maximum over $a \in (0, \infty)$ in $d^*$.   □

For later use we also collect the following properties of $v_{c^*}$ and $\bar{v}_{d^*}$:

LEMMA 3.   Let $x, a > 0$. The following are true:

(i) $v'_{c^*}(x) \geq 1$.

(ii) $1 \leq \bar{v}'_{d^*}(x) \leq \varphi$. Further, if $d^* > 0$, $\bar{v}'_{d^*}(d^{*-}) = 1$ and $\bar{v}'_{d^*}(0^+) = \varphi$ [resp. $\bar{v}'_{d^*}(0^+) < \varphi$] if $X$ has unbounded [resp. bounded] variation.

(iii) $a \mapsto \bar{v}_a(x)$ is nonincreasing for $a > d^*$.

(iv) The function $\bar{v}_a : (0, \infty) \to \mathbf{R} : x \mapsto \bar{v}_a(x)$ is concave.

PROOF.   (i) Since, by Lemma 2, $c^* < \infty$, the statement follows from the definition of $c^*$.

(ii) In view of Lemma 2 and the argument in Proposition 3 it follows that if $d^* > 0$ and $0 < x < d^*$,

$$1 = \varphi Z^{(q)}(x) + W^{(q)}(x)[1 - \varphi Z^{(q)}(x)]/W^{(q)}(x)$$
$$\leq \varphi Z^{(q)}(x) + W^{(q)}(x)[1 - \varphi Z^{(q)}(d^*)]/W^{(q)}(d^*) = \bar{v}'_{d^*}(x).$$

Also, if $d^* > 0$ and $0 < x < d^*$, it holds that

$$(\bar{v}'_{d^*}(x) - \varphi)W^{(q)}(d^*)$$
$$= \varphi(Z^{(q)}(x) - 1)W^{(q)}(d^*) + W^{(q)}(x)[1 - \varphi Z^{(q)}(d^*)]$$
$$= \varphi q[\overline{W}^{(q)}(x)W^{(q)}(d^*) - W^{(q)}(x)\overline{W}^{(q)}(d^*)] + W^{(q)}(x)(1 - \varphi) \leq 0,$$

where in the third line we used Lemma 1. As $W^{(q)}(d^*) > 0$, we conclude that $v'_{d^*}(x) \leq \varphi$. The other statements of (ii) follow from the definitions of $\bar{v}_{d^*}$ and $Z^{(q)}$ and the form of $W^{(q)}(0)$ [see (3.5)].

(iii) The assertion follows since, from the proof of Proposition 3, $(\partial \bar{v}_a(x)/\partial a)$ has the same sign as $F(a)$ and $F(a) \leq 0$ for $a > d^*$.

(iv) Suppose that $d^* > 0$. In view of the definitions of $\overline{Z}^{(q)}$ and $Z^{(q)}$ it holds for $0 < x < d^*$ that

$$\frac{\bar{v}''_{d^*}(x)}{W^{(q)}(x)} = \varphi q + \frac{W^{(q)\prime}(x)}{W^{(q)}(x)}\left[\frac{1 - \varphi Z^{(q)}(d^*)}{W^{(q)}(d^*)}\right]$$
$$\leq \varphi q + \frac{W^{(q)\prime}(d^*)}{W^{(q)}(d^*)}\left[\frac{1 - \varphi Z^{(q)}(d^*)}{W^{(q)}(d^*)}\right] = 0,$$

where we used in the second line that $\varphi Z^{(q)}(d^*) > 1$ and that $W^{(q)\prime}/W^{(q)}$ is decreasing [see the remark just before (3.15)]. Since $\bar{v}'_{d^*}(d^{*-}) = 1$ and $\bar{v}'_{d^*}(x) = 1$ for $x > d^*$, the assertion follows. $\square$

5.4. *Verification theorems.* To investigate the optimality of the barrier strategy $\pi_{c^*}$ across all admissible strategies $\Pi$ for the classical dividend problem (2.2) we are led, by standard Markovian arguments, to consider the following variational inequality:

$$(5.8) \quad \max\{\Gamma w(x) - qw(x), 1 - w'(x)\} = 0, \qquad x > 0,$$
$$w(x) = 0, \qquad x < 0,$$

for functions $w : \mathbf{R} \to \mathbf{R}$ in the domain of the extended generator $\Gamma$ of $X$. For the "bail-out" problem (2.4) we are led to the variational inequality

$$(5.9) \quad \max\{\Gamma w(x) - qw(x), 1 - w'(x)\} = 0, \qquad x > 0,$$
$$w'(x) \leq \varphi, \qquad x > 0, \qquad w'(x) = \varphi, \qquad x < 0,$$

for $w : \mathbf{R} \to \mathbf{R}$ in the domain of $\Gamma$.

To establish the optimality of the barrier strategies among all admissible strategies we shall show the following verification results. In the case of (5.8) we shall only prove a local verification theorem.

PROPOSITION 4. *Let $w : [0, \infty) \to \mathbf{R}$ be continuous.*

(i) *Let $C \in (0, \infty]$, suppose $w(0) = w(0^+) \geq 0$ and extend $w$ to the negative half-line by setting $w(x) = 0$ for $x < 0$. Suppose $w$ is $C^2$ on $(0, C)$ [if $X$ has unbounded variation] or is $C^1$ on $(0, C)$ [if $X$ has bounded variation]. If $w$ satisfies (5.8) for $x \in (0, C)$, then $w \geq \sup_{\pi \in \Pi_{\leq C}} v_\pi$. In particular, if $C = \infty$, $w \geq v_*$.*

(ii) *Suppose $w \in C^2[0, \infty)$ and extend $w$ to the negative half-axis by setting $w(x) = w(0) + \varphi x$ for $x \leq 0$. If $w$ satisfies (5.9), then $w \geq \bar{v}_*$.*

The proof follows below. Inspired by properties of $v_{c^*}$ and with the smoothness required to apply the appropriate version of Itô's lemma in mind, we weaken now the assumptions of the above proposition on the solution $w$. Let $\mathcal{P} = (p_1, p_2, \ldots, p_N)$ with $0 < p_1 < \cdots < p_N$ be a finite subset of $(0, \infty)$ and let $w : [0, \infty) \to [0, \infty)$ be continuous. If $X$ has bounded variation, suppose that $w \in C^1(0, \infty) \backslash \mathcal{P}$ with finite left- and right-derivatives $w'_-(x), w'_+(x)$ for $x \in \mathcal{P}$ and that $w$ satisfies the HJB equation (5.8) where $w'$ is understood to be $w'_-$. If $X$ has unbounded variation, suppose that $w \in C^2(0, \infty) \backslash \mathcal{P}$ with $w(0) = 0$ and finite left- and right-second derivatives for $x \in \mathcal{P}$ and that $w$ satisfies the HJB equation (5.8) where $w''$ is understood to be the weak derivative of $w'$. Finally, we impose a linear growth condition on $w$:

$$(5.10) \qquad\qquad w(x) = \mathcal{O}(x) \qquad (x \to \infty).$$

The following result complements Theorem 9.4 in [1]:

PROPOSITION 5. *Suppose $w$ is as described above. If $w'(0^+) > 1$, then $c^* > 0$ and $w(x) = v_*(x) = v_{c^*}(x)$ for $x \in [0, c^*]$.*

PROOF OF PROPOSITION 4. (ii) Let $\bar{\pi} \in \overline{\Pi}$ be any admissible policy and denote by $L = L^{\bar{\pi}}, R = R^{\bar{\pi}}$ the corresponding pair of cumulative dividend and cumulative loss processes, respectively, and by $V = V^{\bar{\pi}}$ the corresponding risk process. By an application of Itô's lemma to $e^{-qt} w(V_t)$ it can be verified that

$$
\begin{aligned}
(5.11) \quad & e^{-qt} w(V_t) - w(V_0) \\
& = J_t + \int_0^t e^{-qs} w'(V_{s-}) \, dR_s^c - \int_0^t e^{-qs} w'(V_{s-}) \, dL_s^c \\
& \quad + \int_0^t e^{-qs} (\Gamma w - qw)(V_{s-}) \, ds + M_t,
\end{aligned}
$$

where $M_t$ is a local martingale with $M_0 = 0$, $R^c$ and $L^c$ are the pathwise continuous parts of $R$ and $L$, respectively, and $J_t$ is given by

$$(5.12) \qquad J_t = \sum_{s \leq t} e^{-qs} [w(A_s + B_s) - w(A_s)] \mathbf{1}_{\{B_s \neq 0\}},$$

where $A_s = V_{s^-} + \Delta X_s$ and $B_s = \Delta(R - L)_s$ denotes the jump of $R - L$ at time $s$. Note that $1 \leq w'(x) \leq \varphi$ holds for all $x \in \mathbf{R}$. In particular, we see that $w(A_s + B_s) - w(A_s) \leq \varphi \Delta R_s - \Delta L_s$, so that the first three terms on the right-hand side of (5.11) are bounded above by $\varphi \int_0^t e^{-qs} \, dR_s - \int_0^t e^{-qs} \, dL_s$. Let $T_n$ be the first time the absolute value of any of the five terms on the right-hand side of (5.11) exceeds the value $n$, so that, in particular, $T_n$ is a localizing sequence for $M$ with $T_n \to \infty$ a.s. Taking (5.11) at $T_n$, taking expectations and using that, on $[0, \infty)$, $w$ is bounded below by some constant, $-M$ say, $1 \leq w'(x) \leq \varphi$ and $(\Gamma w - qw)(x) \leq 0$ for $x > 0$, it follows after rearranging that

$$w(x) \geq \mathbf{E}_x \left[ \int_0^{T_n} e^{-qs} \, dL_s - \varphi \int_0^{T_n} e^{-qs} \, dR_s \right] + \mathbf{E}_x [e^{-qT_n} w(V_{T_n})]$$

$$\geq \mathbf{E}_x \left[ \int_0^{T_n} e^{-qs} \, dL_s - \varphi \int_0^\infty e^{-qs} \, dR_s \right] - M \mathbf{E}_x [e^{-qT_n}].$$

Letting $n \to \infty$, in view of the fact that $q > 0$ and $T_n \to \infty$ a.s. and the condition (2.3) in conjunction with the monotone convergence theorem, it follows that $\bar{v}_{\bar{\pi}}(x) \leq w(x)$. Since $\bar{\pi}$ was arbitrary we conclude that $w$ dominates the value function $\bar{v}_*$.

(i) Let $\pi \in \Pi_{\leq C}$ and denote by $L = L^\pi$ and $U = U^\pi$ the corresponding cumulative dividend process and risk process, respectively. If $X$ has unbounded variation, $w$ is $C^2$ and we are allowed to apply Itô's lemma (e.g., [23], Theorem 32) to $e^{-q(t \wedge \sigma^\pi)} w(U_{t \wedge \sigma^\pi})$, using that $U_t \leq C$. If $X$ has bounded variation, $w$ is $C^1$ and we apply the change of variable formula (e.g., [23], Theorem 31). Following then an analogous line of reasoning as in (ii) we find that

$$(5.13) \qquad w(x) \geq \mathbf{E}_x \left[ \int_0^{T'_n \wedge \sigma^\pi} e^{-qs} \, dL_s \right]$$

for some increasing sequence of stopping times $T'_n$ with $T'_n \to \infty$ a.s. Taking $n \to \infty$ in (5.13) yields, in view of the monotone convergence theorem and the fact that $w \geq 0$, that,

$$w(x) \geq \mathbf{E}_x \left[ \int_0^{\sigma^\pi} e^{-qs} \, dL_s \right].$$

Since the previous display holds for arbitrary $\pi \in \Pi_{\leq C}$, it follows that $w(x) \geq \sup_{\pi \in \Pi_{\leq C}} v_\pi(x)$ and the proof is complete. $\square$

PROOF OF PROPOSITION 5. Noting that $w$ is smooth enough for an application of the appropriate version of Itô's lemma (as follows from the (proof) of the Itô–Tanaka–Meyer formula; see, e.g., [23]), it can be verified, as in Proposition 4, that $w \geq v_*$.

Putting $m = \inf\{x > 0 : w'(x^-) = 1\}$, it follows from the assumptions that $m \in (0, \infty)$ or $m = \infty$. The latter case can be ruled out as follows. In view of the facts that $w(0^+) > 1$ and $w$ satisfies (5.8), the assumption that $m = \infty$ implies that $w'(x) > 1$ and $\Gamma w(x) - q w(x) = 0$ for $x > 0$. An application of Itô's lemma to $e^{-q(t \wedge T_{0,a})} w(X_{t \wedge T_{0,a}})$ shows then that for $0 < x < a$ it holds that

$$w(x) = \mathbf{E}_x[e^{-q(T_{0,a})} w(X_{T_{0,a}})] = W^{(q)}(x)[w(a)/W^{(q)}(a)],$$

where we used (3.6) and that $w(x) = 0$ for $x < 0$ and $w(0) = 0$ if $\sigma > 0$ (as $\mathbf{P}_x[X_{T_0^-} = 0] > 0$ if $\sigma > 0$). Letting $a \to \infty$, the right-hand side converges to zero in view of the linear growth condition (5.10) and the fact that $W^{(q)}(a)$ grows exponentially fast as $a \to \infty$. Since $w \geq v_*$, this leads to a contradiction, and we see that $m \in (0, \infty)$. Applying subsequently Itô's lemma to $e^{-q(t \wedge \sigma^\pi)} w(U_{t \wedge \sigma^\pi})$ with $\pi = \pi_m$ and using that $w'(x) > 1$ and $\Gamma w(x) - q w(x) = 0$ for $x \in (0, m)$, we find that

$$w(x) = \mathbf{E}_x\left[ \int_0^{T_n'' \wedge \sigma^{\pi_m}} e^{-qs} \, dL_s \right] + \mathbf{E}_x[e^{-q(\sigma^{\pi_m} \wedge T_n'')} w(U_{\sigma^{\pi_m} \wedge T_n''})]$$

for some increasing sequence of stopping time $T_n''$ with $T_n'' \to \infty$. Letting $n \to \infty$ and using that $w(U_{\sigma^{\pi_m} \wedge T_n''})$ is bounded (since $U^{\pi_m} \leq m$) and $w(U_{\sigma^{\pi_m}}) = 0$, it follows that $w(x) = v_m(x)$ for $x \in [0, m]$. Since, on the one hand, Proposition 3 implies that $v_m \leq v_{c^*}$, while, on the other hand, $w \geq v_*$, we deduce that $m = c^*$ and $v_*(x) = v_{c^*}(x)$ for $x \in [0, c^*]$ with $c^* > 0$.   □

5.5. *Proofs of Theorems 2 and 3.* We set $v_{c^*}(x) = 0$ for $x < 0$ and extend $\bar{v}_{d^*}$ to the negative half-axis by setting $\bar{v}_{d^*}(x) = \bar{v}_{d^*}(x) + \varphi x$ for $x < 0$. Recalling that $W^{(q)}(x) = 0$, $Z^{(q)}(x) = 1$ and $\overline{Z}^{(q)}(x) = x$ for $x < 0$, we see that these are natural extensions of the formulas (5.1) and (5.4) and satisfy the HJB equations (5.8) and (5.9) for $x < 0$. The proofs of Theorems 2 and 3 are based on the following lemmas:

LEMMA 4. *If $c^* > 0$, $(\Gamma v_{c^*} - q v_{c^*})(x) = 0$ for $x \in (0, c^*)$.*

LEMMA 5. *It holds that $(\Gamma \bar{v}_{d^*} - q \bar{v}_{d^*})(x) \leq 0$ [resp. $= 0$] if $x > 0$ [resp. if $d^* > 0$ and $x \in (0, d^*)$].*

PROOF OF THEOREM 2. (i) In view of Lemmas 3(i) and 4 it follows that the function $v_{c^*}$ satisfies the variational inequality (5.8) for $x \in (0, c^*)$.

Therefore, Proposition 4 implies the optimality of the strategies $\pi_c^*$ in the set $\Pi_{\leq c^*}$.

(ii) If the condition of Theorem 2(ii) holds, then, in view of Lemmas 3(i) and 4, it follows that $v_{c^*}$ satisfies the variational inequalities (5.8) for $x \in (0, \infty)$. By Proposition 4 it then follows that $v_c^* = v_*$. □

PROOF OF THEOREM 3. In view of Lemmas 3(ii) and 5 it follows that the function $\bar{v}_{d^*}$ satisfies the variational inequality (5.9). Therefore, Proposition 4(ii) implies that $\bar{v}_* = \bar{v}_{d^*}$ and the strategy $\bar{\pi}_{0,d^*}$ is optimal. □

PROOF OF LEMMA 4. Suppose that $c^* > 0$. Since $e^{-q(t \wedge T_{0,c^*})} W^{(q)}(X_{t \wedge T_{0,c^*}})$ is a martingale, $e^{-q(t \wedge T_{0,c^*})} v_{c^*}(X_{t \wedge T_{0,c^*}})$ inherits this martingale property by definition of $v_{c^*}$. Since $v_{c^*}$ is smooth enough to apply the appropriate version of Itô's lemma ([23], Theorem 31) is applicable if $X$ has bounded variation since then $v_{c^*} \in C^1(0, c^*)$ and [23], Theorem 32 is applicable if $X$ has unbounded variation as then $v_{c^*} \in C^2(0, c^*)$, the martingale property in conjunction with Itô's lemma implies that $\Gamma v_{c^*}(y) - q v_{c^*}(y) = 0$ for $0 < y < c^*$. □

PROOF OF LEMMA 5. First let $d^* > 0$. In view of Proposition 2 and the martingale property (3.7), it follows that $e^{-q(t \wedge T_{0,d^*})} \bar{v}_{d^*}(X_{t \wedge T_{0,d^*}})$ is a martingale. An application of Itô's lemma, which we are allowed to apply as $Z^{(q)} \in C^2(0, \infty)$, then yields that $\Gamma \bar{v}_{d^*}(y) - q \bar{v}_{d^*}(y) = 0$ for $0 < y < d^*$.

Let now $d^* \geq 0$. From the form of the infinitesimal generator $\Gamma$ and the definition of $\bar{v}_{d^*}$ it follows that for $x > d^*$ we have

$$\Gamma \bar{v}_{d^*}(x) - q \bar{v}_{d^*}(x)$$
$$= c \bar{v}_{d^*}'(x) + \int_{-\infty}^0 \{\bar{v}_{d^*}(x+y) - \bar{v}_{d^*}(x) - y\mathbf{1}_{\{|y|<1\}} \bar{v}_{d^*}'(x)\} \nu(dy) - q \bar{v}_{d^*}(x)$$
$$= c + \int_{-\infty}^0 \{\bar{v}_{d^*}(x+y) - (x-b) - y\mathbf{1}_{\{|y|<1\}}\} \nu(dy) - q(x-b)$$

where $c$ is some constant and $b = \bar{v}_{d^*}(d^*) - d^*$. Since $\bar{v}_{d^*}$, considered as function mapping $\mathbf{R}$ to $\mathbf{R}$, is concave [Lemma 3(ii)–(iv)] and integration preserves concavity, it follows from the previous display that $g|_{(d^*, \infty)}$ with $g : (0, \infty) \to \mathbf{R}$ given by $g(x) = \Gamma \bar{v}_{d^*}(x) - q \bar{v}_{d^*}(x)$ is concave. Further, in view of continuity of $g$ and the fact that $g|_{(0,d^*)} = 0$ it follows that $\lim_{x \downarrow d^*} g(x) = \lim_{x \uparrow d^*} g(x) = 0$. We claim that the right-derivative of $g$ in $d^*$ is nonpositive, $g_+'(d^*) \leq 0$. Before we prove this claim we first show that the claim implies that $g(x) \leq 0$ for $x > d^*$. Indeed, if $g_+'(d^*) \leq 0$ then, in view of the concavity of $g$, it holds that $g_+'(x) \leq 0$ for $x > d^*$ so that $g(x) \leq g(d^*) = 0$ for $x > d^*$.

The rest of the proof is devoted to showing that $g'_+(d^*) \le 0$. To that end, we shall first show that, for $a > d^*$ and $V_0 = x \in (0, a)$, it holds that

$$
\begin{aligned}
\bar{v}_a(x) - \bar{v}_{d^*}(x) &= \mathbf{E}_x\left[\int_0^\infty e^{-qs}(\Gamma\bar{v}_{d^*} - q\bar{v}_{d^*})(V_{s^-}^a)\,ds\right] \\
&= \int_{d^*}^a (\Gamma\bar{v}_{d^*} - q\bar{v}_{d^*})(y)\widetilde{U}_a^q(x, dy),
\end{aligned}
\tag{5.14}
$$

where $\widetilde{U}_a^q(x, dy) = \widetilde{U}^q(x, dy)$ is the resolvent measure of $V^a$ given in (4.1).

Note that $s \mapsto L_s^a$ can be taken to be continuous in this case and that the support of Stieltjes measure $dL_s^a$ is contained in the set $\{s : V_{s^-}^a = a\}$. Further, in this case $R^0$ jumps at time $s$ if and only if $X$ jumps at time $s$ and $\Delta X_s$ is larger than $V_{s^-}^a$. Thus $\Delta R_s^0 = -\min\{0, V_{s^-}^a + \Delta X_s\}$ and the measure $d(R^0)_s^c$ has support inside $\{s : V_{s^-}^a = 0\}$. In view of these observations, an application of Itô's lemma to $e^{-qt}\bar{v}_{d^*}(V_t^a)$ as in (5.11) shows that

$$
\begin{aligned}
e^{-qt}\bar{v}_{d^*}(V_t^a) &- \bar{v}_{d^*}(x) \\
&= \int_0^t e^{-qs}\bar{v}'_{d^*}(0^+)\,d(R^0)_s^c \\
&\quad - \int_0^t e^{-qs}\bar{v}'_{d^*}(a^-)\,dL_s^a \\
&\quad + \varphi\sum_{s \le t} e^{-qs}\Delta R_s^0\mathbf{1}_{\{\Delta R_s^0 > 0\}} \\
&\quad + \int_0^t e^{-qs}(\Gamma\bar{v}_{d^*} - q\bar{v}_{d^*})(V_{s^-}^a)\,ds + M_t,
\end{aligned}
\tag{5.15}
$$

where we used that in (5.11) $J_t = \varphi\sum_{s \le t} e^{-qs}\Delta R_s^0\mathbf{1}_{\{\Delta R_s^0 > 0\}}$ since, by definition of the extended function $\bar{v}_{d^*}$ on $(-\infty, 0]$, it follows that $\bar{v}_{d^*}(x + y) - \bar{v}_{d^*}(x) = \varphi y$ if $x = -y$, $x < 0$. Since $\bar{v}_{d^*} \in C^2(0, \infty)$ and $V^a$ takes values in $[0, a]$, it follows that in this case $M$ is a martingale. Further, it holds that $\bar{v}'_{d^*}(a^-) = 1$ (recalling that $a > d^*$) and either $(R^0)^c = 0$ (if $X$ has bounded variation) or $\bar{v}'_{d^*}(0^+) = \varphi$ (if $X$ has unbounded variation). Taking then expectations and letting $t \to \infty$ in (5.15) shows, in view of the dominated convergence theorem and the first part of the proof, that (5.14) holds true.

We now finish the proof by supposing $g'_+(d^*) > 0$ and showing that this assumption leads to a contradiction. If $g'_+(d^*) > 0$, then, in view of the continuity of $g$, there exists an $\varepsilon > 0$ such that $g(x) > 0$ for $x \in (d^*, d^* + \varepsilon)$. Since, for $a > 0$, $\widetilde{U}_a^q(x, dy)$ is absolutely continuous on $(0, a)$ with positive density [see (4.1)], it follows from (5.14) that $\bar{v}_{d^* + \varepsilon}(x) > \bar{v}_{d^*}(x)$, which contradicts Proposition 3(ii). Therefore $g'_+(d^*) \le 0$ and the proof is done. □

## 6. Examples.

6.1. *Small claims*: *Brownian motion*. If $X_t = \sigma B_t + \mu t$ is a Brownian motion with drift $\mu$ (a standard model for small claims), then

$$W^{(q)}(x) = \frac{1}{\sigma^2 \delta}[e^{(-\omega+\delta)x} - e^{-(\omega+\delta)x}],$$

where $\delta = \sigma^{-2}\sqrt{\mu^2 + 2q\sigma^2}$ and $\omega = \mu/\sigma^2$. It is a matter of calculus to verify that

$$W^{(q)\prime\prime}(x) = 2\sigma^{-2}[qW^{(q)}(x) - \mu W^{(q)\prime}(x)]$$

from which it follows that if $\mu \leq 0$, $W^{(q)\prime}(x)$ attains its minimum over $[0, \infty)$ in $x = 0$. Thus in the classical setting it is optimal to take out all dividends immediately if $\mu \leq 0$; if $\mu > 0$ it follows that $c^* > 0$ and it holds that $W^{(q)\prime\prime}(c^*) = 0$, so that $W^{(q)}(c^*)/W^{(q)\prime}(c^*) = \mu/q$, as Gerber and Shiu [12] have found before, and the optimal level $c^*$ is explicitly given by

$$c^* = \log\left|\frac{\delta + \omega}{\delta - \omega}\right|^{1/\delta}.$$

Since $\frac{\sigma^2}{2}v''_{c^*}(x) + \mu v'_{c^*}(x) - qv_{c^*}(x) < 0$ for $x > c^*$, it follows by Theorem [2] that $\pi_{c^*}$ is the optimal strategy as shown before by Jeanblanc and Shiryaev [15]. In the "bail-out" setting $d^* \in (0, \infty)$ solves $G(a) = 0$ where $G$ is given in (5.6) with

$$Z^{(q)}(y) = \frac{q}{\sigma^2 \delta}\left[\frac{1}{\omega + \delta}e^{-(\omega+\delta)y} + \frac{1}{\delta - \omega}e^{(-\omega+\delta)y}\right]$$

and

$$W^{(q)\prime}(y) = \frac{1}{\sigma^2 \delta}[(\omega + \delta)e^{-(\omega+\delta)y} + (\delta - \omega)e^{(-\omega+\delta)y}].$$

The relation between the classical and bail-out strategies in this Brownian setting is studied in [19].

6.2. *Stable claims.* We model $X$ as

$$X_t = \sigma Z_t,$$

where $Z$ is a standard stable process of index $\alpha \in (1, 2)$ and $\sigma > 0$. Its cumulant is given by $\psi(\theta) = (\sigma\theta)^\alpha$. By inverting the Laplace transform $(\psi(\theta) - q)^{-1}$, Bertoin [5] found that the $q$-scale function is given by

$$W^{(q)}(y) = \alpha\frac{y^{\alpha-1}}{\sigma^\alpha}E'_\alpha\left(q\frac{y^\alpha}{\sigma^\alpha}\right), \qquad y > 0,$$

and hence $Z^{(q)}(y) = E_\alpha(q(y/\sigma)^\alpha)$ for $y > 0$, where $E_\alpha$ is the Mittag–Leffler function of index $\alpha$

$$E_\alpha(y) = \sum_{n=0}^{\infty} \frac{y^n}{\Gamma(1 + \alpha n)}, \qquad y \in \mathbf{R}.$$

The form of the value functions $v_{c^*}$ and $\bar{v}_{d^*}$ follows by inserting the expressions for the scale functions in (5.1)–(5.4). The optimal levels $c^*, d^*$ are given by

$$c^* = \sigma q^{-1/\alpha} u(\alpha)^{1/\alpha}, \qquad d^* = \sigma q^{-1/\alpha} v(\alpha)^{1/\alpha},$$

where $u(\alpha) > 0$ and $v(\alpha) > 0$ are positive roots of the respective equations

$$(\alpha - 1)(\alpha - 2)E'_\alpha(u) + 3\alpha(\alpha - 1)uE''_\alpha(u) + \alpha^2 u^2 E'''_\alpha(u) = 0,$$

$$\varphi\alpha v(E'_\alpha(v))^2 + [(\alpha - 1)E'_\alpha(v) + \alpha v E''_\alpha(v)][1 - \varphi E_\alpha(v)] = 0.$$

6.3. *Cramér–Lundberg model with exponential jumps.* Suppose $X$ is given by the Cramér–Lundberg model (1.1) with exponential jump sizes, that is, $X$ is a deterministic drift $p$ (the premium income) minus a compound Poisson process (with jump intensity $\lambda$ and jump sizes $C_k$ that are exponentially distributed with mean $1/\mu$) such that $X$ has positive drift; that is, $p > \lambda/\mu$. Then $\psi(\theta) = p\theta - \lambda\theta/(\mu + \theta)$ and the scale function $W^{(q)}$ is given by

$$W^{(q)}(x) = p^{-1}(A_+ e^{q^+(q)x} - A_- e^{q^-(q)x}),$$

where $A_\pm = \frac{\mu + q^\pm(q)}{q^+(q) - q^-(q)}$ with $q^+(q) = \Phi(q)$ and $q^-(q)$ the smallest root of $\kappa(\theta) = q$:

$$q^\pm(q) = \frac{q + \lambda - \mu p \pm \sqrt{(q + \lambda - \mu p)^2 + 4pq\mu}}{2p}.$$

Then from (5.3) we have that $c^* = 0$ if $W^{(q)''}(0) \geq 0 \Leftrightarrow p\lambda\mu \leq (q + \lambda)^2$. If $p\lambda\mu > (q + \lambda)^2$, $c^*$ is strictly positive and given by

$$c^* = \frac{1}{q^+(q) - q^-(q)} \log \frac{q^-(q)^2(\mu + q^-(q))}{q^+(q)^2(\mu + q^+(q))}.$$

Since it is readily verified that $\Gamma v_{c^*}(x) - q v_{c^*}(x) < 0$ for $x > c^*$, Theorem 2(ii) implies that $\pi_{c^*}$ is the optimal strategy.

Further, if $\lambda(\varphi - 1) \leq q$, then $d^* = 0$. Otherwise $d^* > 0$ satisfies $G(d^*) = 0$ where $G$ is given in (5.6).

6.4. *Jump-diffusion with hyper exponential jumps.* Let $X = \{X_t, t \geq 0\}$ be a jump-diffusion given by

$$X_t = \mu t + \sigma W_t - \sum_{i=1}^{N_t} Y_i,$$

where $\sigma > 0$, $N$ is a Poisson process with intensity $\lambda > 0$ and $\{Y_i\}$ is a sequence of i.i.d. random variables with hyper exponential distribution

$$F(y) = 1 - \sum_{i=1}^{n} A_i e^{-\alpha_i y}, \qquad y \geq 0,$$

where $A_i > 0$; $\sum_{i=1}^{n} A_i = 1$; and $0 < \alpha_1 < \cdots < \alpha_n$. In [4] it was shown that the function $Z^{(q)}$ of $X$ is given by

$$Z^{(q)}(x) = \sum_{i=0}^{n+1} D_i(q) e^{\theta_i(q)x},$$

where $\theta_i = \theta_i(q)$ are the roots of $\psi(\theta) = q$, where $\theta_{n+1} > 0$ and the rest of the roots are negative, and where

$$D_i(q) = \prod_{k=1}^{n} (\theta_i(q)/\alpha_k + 1) \Big/ \prod_{k=0, k \neq i}^{n+1} (\theta_i(q)/\theta_k(q) - 1).$$

If $c^* > 0$, it is a nonnegative root $x$ of

$$\sum_{i=0}^{n+1} \theta_i(q)^3 D_i(q) e^{\theta_i(q)x} = 0.$$

**Acknowledgments.** We would like to thank the anonymous referees for their careful review and valuable comments. We would like to thank Hansjörg Albrecht for his comments that led to substantial improvements of the paper.

F. AVRAM                                   Z. PALMOWSKI
DEPARTEMENT DE MATHEMATIQUES               DEPARTMENT OF MATHEMATICS
UNIVERSITÉ DE PAU                          UNIVERSITY OF WROCŁAW
64000 PAU                                  PL. GRUNWALDZKI 2/4
FRANCE                                     50-384 WROCŁAW
E-MAIL: Florin.Avram@univ-Pau.fr           POLAND
                                           E-MAIL: zpalma@math.uni.wroc.pl

M. R. Pistorius
Department of Mathematics
King's College London
The Strand
London WC2R 2LS
United Kingdom
E-mail: Martijn.Pistorius@kcl.ac.uk